\documentclass[runningheads,a4paper]{llncs}

\usepackage{amssymb}
\usepackage{graphicx}
\usepackage{color}

\newcommand{\keywords}[1]{\par\addvspace\baselineskip
\noindent\keywordname\enspace\ignorespaces#1}

\newcommand{\R}{\mathbb R}

\newcommand{\eqref }[1]{(#1)}
\newcommand{\approxT}{\widehat T}

\begin{document}

\mainmatter  

\title{Two Semi--Lagrangian Fast Methods for Hamilton--Jacobi--Bellman Equations\thanks{This research was supported by the following grants: AFOSR Grant FA9550-10-1-0029,  ITN-Marie Curie Grant 264735-SADCO.}}
%
\titlerunning{Two Semi-Lagrangian Fast Methods for HJB Equations}
%
%
%
\author{Simone Cacace$^\dag$ \and Emiliano Cristiani$^\ddag$ \and Maurizio Falcone$^\dag$}
\authorrunning{S. Cacace, E. Cristiani, M. Falcone}
%
\institute{
$^\ddag$Dipartimento di Matematica, Sapienza -- Universit\`a di Roma, Rome, Italy \\
$^\dag$Istituto per le Applicazioni del Calcolo, CNR, Rome, Italy \\
}
%
%
%
%
\toctitle{Lecture Notes in Computer Science}
\tocauthor{Authors' Instructions}
\maketitle
\begin{abstract}
In this paper we apply the Fast Iterative Method (FIM) for solving general Hamilton--Jacobi--Bellman (HJB) 
equations and we compare the results with an accelerated version of the Fast Sweeping Method (FSM). 
We find that FIM can be indeed used to solve HJB equations with no relevant modifications with respect to the original algorithm proposed for the eikonal equation, and that it overcomes FSM in many cases. 
Observing the evolution of the active list of nodes for FIM, we recover another numerical validation of the arguments recently discussed in \cite{CCF14} about the 
impossibility of creating local single-pass methods for HJB equations.
\end{abstract}
\section{Introduction}
The study of Hamilton--Jacobi (HJ) equations arises in several contexts, including classical mechanics, front propagation, control problems and differential games. 
In particular, for optimal control problems, the value function can be characterized as the unique viscosity solution of a Hamilton--Jacobi--Bellman (HJB) equation. 
Unfortunately, solving numerically the HJB equation can be rather expensive from the computational point of view. This is the reason why in the last years an increasing number of efficient techniques have been proposed, see, e.g., \cite{CCF14} for a brief review. 

Basically, these algorithms are divided in two main classes: \emph{single-pass} and \emph{iterative}. An algorithm is said to be \emph{single-pass} if one can fix \textit{a priori} a (small) number $r$ which depends only on the equation and on the mesh structure (not on the number of mesh points) such that each mesh point is re-computed at most $r$ times.
Single-pass algorithms usually divide the numerical grid in, at least, three changing subsets: \emph{Accepted} (ACC), \emph{Considered} (CONS), and \emph{Far} (FAR). 
Nodes in ACC are definitively computed, nodes in CONS are computed but their values are not yet final, and nodes in FAR are not yet computed.
We say that a single-pass algorithm is \emph{local} if the computation at any mesh point involves only the values of first neighboring nodes, the region CONS is 1-cell-thick  and no information coming from FAR region is used.
The methods which are not single-pass are iterative.

Among fast methods, the prototype algorithm for the local single-pass class is the Fast Marching  Method (FMM) \cite{S96,T95}, while that for the iterative class is the Fast Sweeping Method (FSM) \cite{KOQ04,QZZ07,TCOZ04,Z05}. Another interesting method is the Fast Iterative Method (FIM) \cite{FJPKW11,JW08,JW07}, which shares some features with both iterative and single-pass methods.
Recently, Cacace et al.\ \cite{CCF14} have shown that it is not possible to create a local single-pass algorithm for solving general HJB equations. 
This motivates the efforts to develop new techniques, particularly in the class of iterative methods.

In this paper, we focus on the following minimum time HJB equation
\begin{equation}\label{HJB}
\sup\limits_{a\in B_1}\left\{-f(x,a)\cdot\nabla T(x)\right\} = 1, \qquad 
x\in\R^d\backslash\mathcal T
\end{equation}
where $d$ is the space dimension, $\mathcal T$ is a closed nonempty target set in $\R^d$, \mbox{$f:\R^d\times B_1\to\R^d$} is a given vector-valued Lipschitz continuous function, and $B_1$ is the unit ball in $\R^d$ centered in the origin, representing the set of the admissible controls. 
We complement the equation with homogeneous Dirichlet condition $T=0$ on $\mathcal T$.
Let us note that if $f(x,a)=c(x)a$, equation (\ref{HJB}) becomes the eikonal equation 
$
c(x)|\nabla T(x)|=1
$.
To simplify the notations, we restrict the discussion to the case $d=2$. 
Generalizations of the considered algorithms to any space dimension is straightforward, although the implementation is not trivial. 

The goal of this paper is twofold: First, we investigate the possibility of applying a semi-Lagrangian version of the FIM to equation (\ref{HJB}). 
To our knowledge, FIM was only used for solving the eikonal equation \cite{JW08} and a special class of HJ equations \cite{JW07}, 
although there is no particular constraint to apply it in a more general framework. The algorithm indeed does not rely on the special form and features of the eikonal equation. 
In addition, we measure the degree of ``iterativeness'' of FIM, keeping track of how many times each grid node is inserted into the list of nodes which are actually computed at each step. Interestingly, the results indirectly confirm the findings of \cite{CCF14}, showing that general HJB equations require the nodes to be visited an (\textit{a priori}) unknown number of times, i.e.\ single-pass methods do not apply.

Second, we propose a new acceleration technique for the FSM, which is effective when (\ref{HJB}) is discretized by means of a 
semi-Lagrangian scheme (see \cite{FFbook} for a comprehensive introduction). 
It reduces the CPU load for the sup search in (\ref{HJB}), neglecting the control directions which are downwind with respect to the current sweep. 
The new method results to be remarkably faster, although (in general) the number of iterations needed for convergence increases.
\section{Semi-Lagrangian approximation}
Let us  introduce a structured grid $G$ and denote its nodes by $x_i$, $i=1,\ldots,N$. 
The space step is assumed to be uniform and equal to $\Delta x>0$. 
Standard arguments \cite{FFbook} lead to the following discrete version of equation (\ref{HJB}): 
\begin{equation}\label{SLscheme}
T(x_i)\approx \approxT(x_i) = 
\min\limits_{a\in B_1}\left\{ \approxT(\tilde x_{i,a})+ \frac{|x_i-\tilde x_{i,a}|}{|f(x_i,a)|}\right\}\,,\qquad x_i\in G
\end{equation}
where $\tilde x_{i,a}$ is a \textit{non-mesh} point, obtained by integrating, until a certain final time $\tau$, the ordinary differential equation
\begin{equation}\label{ODE}
\left\{
\begin{array}{ll}
\dot y(t)=f(y,a), & \quad  t\in[0,\tau] \\
y(0)=x_i & 
\end{array} 
\right .
\end{equation}
and then setting $\tilde x_{i,a}=y(\tau)$. To make the scheme fully discrete, the set of admissible controls $B_1$ 
is discretized with $N_c$ points and we denote by $a^*$ the optimal control achieving the minimum.
Note that we can get different versions of the semi-Lagrangian (SL) scheme (\ref{SLscheme})
varying $\tau$, the method used to solve (\ref{ODE}), and the interpolation method used to compute  $\approxT (\tilde x_{i,a})$. 
Moreover, we remark that, in any single-pass method, the computation of $\approxT (x_i)$ cannot involve the value $\approxT (x_i)$ itself, because this self-dependency would make the method iterative.
Here we use a 3-point scheme: Equation (\ref{ODE}) is solved by an explicit forward Euler scheme until the solution is at distance 
$\Delta x$ from $x_i$, where it falls inside the triangle of vertices $x_{i,1}$, $x_{i,2}$, and $x_{i,3}$, to be chosen among the first neighbors of 
$x_i$. 
The value $\approxT (x_i)$ is computed by a two-dimensional linear interpolation of the values $\approxT (x_{i,1})$, $\approxT (x_{i,2})$ and 
$\approxT (x_{i,3})$ (see \cite{CCF14} for details). 

\section{Limits of local single-pass methods}\label{sec:limitations}
In this section we briefly recall the main result of \cite{CCF14}. From the numerical point of view, 
it is meaningful to divide HJB equations into four classes. For any given mesh, we have:
\begin{itemize}
\item[(ISO)~] Equations whose characteristic curves coincide or lie in the same simplex of the gradient curves of their solutions.
\item[($\neg$ISO)] Equations for which there exists at least a grid node where the characteristic curve and the gradient curve of the solution do not lie in the same simplex.
\end{itemize}
\begin{itemize}
\item[(REG)~] Equations with non-crossing (regular) characteristic curves. Characteristics spread from the target $\mathcal T$ to the rest of the domain without intersecting.
\item[($\neg$REG)] Equations with crossing characteristic curves. Characteristics start from the target $\mathcal T$ and then meet in finite time, creating shocks. As a result, the solution $T$ is not differentiable at shocks.
\end{itemize}

Let us summarize here the main remarks on single-pass methods: \\
\noindent $\bullet$ FMM works for equations of type ISO and fails for equations of type $\neg$ISO (see \cite{SV03} for further details and explanations), while FSM can be successfully applied in any case.

\noindent $\bullet$ Handling $\neg$ISO case requires CONS not to follow the level sets of the solution itself. Indeed, if CONS turns out to be an approximation of the level sets of the solution, it means that the solution is computed in an increasing order, thus following the gradient curve rather than the characteristic curve.

\noindent $\bullet$ Handling $\neg$REG case requires CONS to be an approximation of the level sets of the solution. Let us clarify this point. Consider the $\neg$REG case and let $x$ be a point belonging to a shock, i.e.\ where the solution is not differentiable. By definition, the value $T(x)$ is carried by two or more characteristic curves reaching $x$ at the same time. Similarly, let $x_i$ be a grid node $\Delta x$-close to the shock. In order to mimic the continuous case,
$x_i$ has to be approached by the ACC region approximately at the same time from the directions corresponding to the characteristic curves. 
In this case, the value $T(x_i)$ is correct (no matter which upwind direction is chosen) and, more important, the characteristic information stops at $x_i$ and it is no longer propagated, getting stuck by the ACC region. As a consequence, the shock is localized properly.
On the other hand, 
if CONS region  is not an approximation of a level set of the solution a node $x_i$ close to a shock can be reached by ACC at different times. When ACC reaches $x_i$ for the first time, it is impossible to detect the presence of the shock by using only local information. 
Indeed, only a global view of the solution allows one to know that another characteristic curve will reach $x_i$ at a later time. 
As a consequence, the algorithm continues the enlargement of CONS and ACC, thus making an error that cannot be corrected by the following 
iterations.

In conclusion, we get that local single-pass methods cannot handle equations $\neg$ISO \& $\neg$REG. In this situation, one has to add non local information regarding the location of the shock, or going back to nodes in ACC at later time, breaking the single-pass property.
This motivates the investigation of new techniques, especially iterative methods, as the ones described in the next sections.

\section{Fast iterative method}
In this section we briefly recall the construction of FIM \cite{FJPKW11,JW08,JW07}.
As in FMM, the main idea of FIM is to update only few grid nodes at each step. These nodes are stored in a separated list, called \emph{active} list. 
During each step, the list of active nodes is modified, and the band
thickens or expands to include all nodes that could be affected by the current updates.
A node can be removed from the active list when its value is up to date with respect to its neighbors (i.e., it has reached convergence) 
and can be appended to the list (relisted) whenever any upwind neighbor's value has changed. 

FIM is formally an iterative method, since the number of times a grid node is visited depends on the dynamics and on the grid size. On the other hand, the active list resembles the set CONS of FMM, and in some special cases FIM is in fact a single-pass algorithm, see Section \ref{sec:tests}. Nevertheless, the active list and CONS 
differ for some important features. The first is that the active list is not kept ordered, 
and then the causality relationship among grid nodes is lost. The second is that the active list can be more than 1-node-thick, i.e. it can approximate a two-dimensional set. Finally, grid nodes removed from the active list can re-enter at a later time. This is the price to pay for loosing the causality.

The FIM algorithm consists of two parts, the initialization and the updating. 
In the initialization step, one has to set the boundary conditions and set the values of the rest of the grid nodes to infinity (or some very large value). Next, the adjacent neighbors of the source nodes (i.e.\ the target) are added to the active list. 
In the updating step, for every point in the list, one computes the new value and checks if the value at the node has converged by comparing the old and the new value at the considered point. If it has converged, one removes the node from the list and append to the list any non active adjacent node such that its updated value is less than the current one. 
The algorithm runs until the list is empty. 

FIM was introduced for solving a special class of HJ equations \cite{JW08,JW07}. 
Nevertheless, in Section \ref{sec:tests} we show that FIM based on a SL discretization \textit{can be successfully applied to general HJB equations with no modifications}. 
\section{An optimized fast sweeping method}
FSM is another popular method for solving HJ equations \cite{KOQ04,QZZ07,TCOZ04,Z05}. 
The main advantage of the method is its implementation, which is extremely easy (easier than that of FMM and FIM). 
FSM is basically the classical iterative (fixed-point) method, since each node is visited in a predefined order, until convergence is reached. 
Here, the visiting directions (sweeps) are alternated in order to follow all possible characteristic directions, trying to exploit causality. 
In two-dimensional problems, the grid is visited sweeping in four directions:  
$S\to N$ \& $W\to E$,
$S\to N$ \& $E\to W$, 
$N\to S$ \& $E\to W$ and
$N\to S$ \& $W\to E$. \\ 
The key point is the Gauss-Seidel-like update of grid nodes, which allows one to compute in a cascade fashion a relevant part of the 
grid nodes in only one sweep. Indeed it is well known that in the case of eikonal equations FSM converges in only four sweeps \cite{Z05}.

\begin{figure}[!h]
 \begin{center}
  \includegraphics[width=0.8\textwidth]{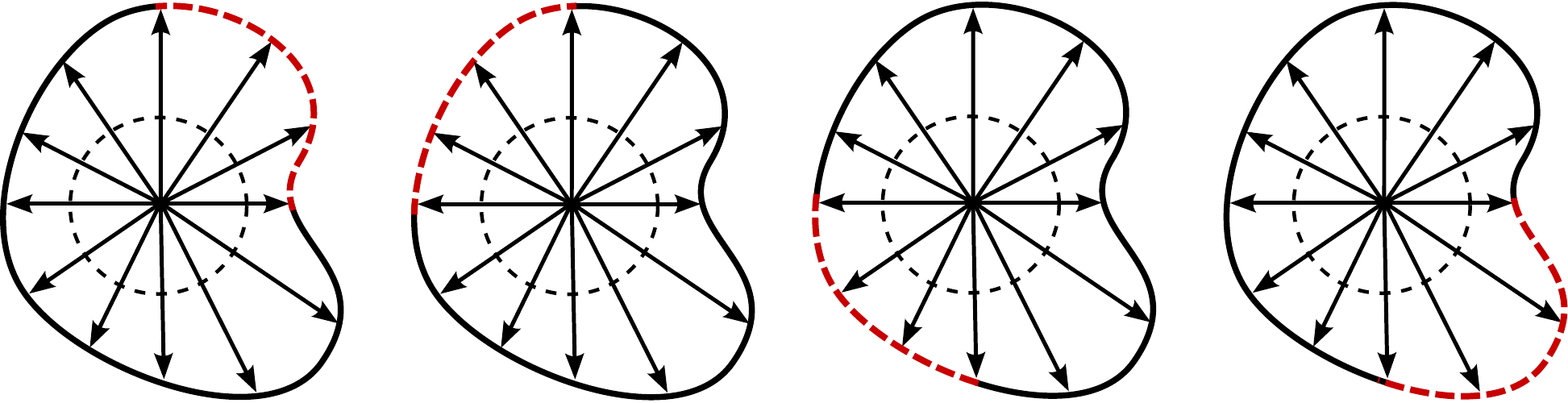}
 \end{center}
 \vskip-.5cm
\caption{Downwind controls with respect to the four sweeps: dashed arcs identify the directions to be removed}
\label{fig:upstreamcontrols}
\end{figure}

Here we propose an easy modification of the FSM based on a SL discretization, aiming at saving CPU time for each sweep. 
Let us explain the idea in the case of a dynamics of the form $c(x,a)a$, with $c>0$. 
It is clear that during the sweep $S\to N$ \& $W\to E$ the algorithm cannot exploit the power of the Gauss-Seidel cascade for the information coming 
from NE. Indeed, even if a node actually depends on its NE neighbor, that information flows upwind and it is not propagated to other nodes during 
the current sweep. Then, we propose to \textit{remove downwind discrete controls from the minimum search in the SL scheme} (\ref{SLscheme}), 
since they have small or no effect in the update of the nodes, see Figure \ref{fig:upstreamcontrols}.
\\
The assumption $c>0$ is needed to preserve the order of the quadrants between the control $a$ and the resulting dynamics $c(x,a)a$. 
Otherwise, the choice of controls to be removed should be adapted according to the sign of $c$.\\
Note that the control set $B_1$ is reduced to a upwind $3/4$ of ball. 
Then, let us denote by Upwind Fast Sweeping Method 3/4 (UFSM{\footnotesize 3/4}) the classical FSM with this control reduction.
An additional speedup, that we expect to work only in case where the characteristics are essentially straight, consists in 
reducing further the control ball to a upwind $1/4$ of ball (UFSM{\footnotesize 1/4}).
\section{Numerical experiments}\label{sec:tests}
\begin{figure}[!h]
\centering
 \begin{tabular}{ccc}
  \includegraphics[width=.3\textwidth]{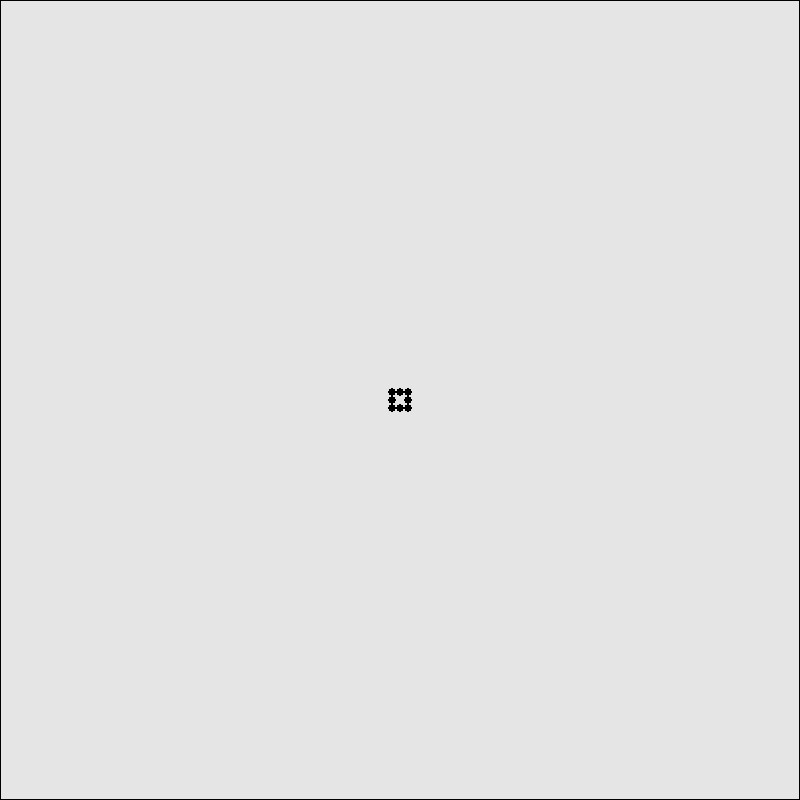} \qquad & 
  \includegraphics[width=.3\textwidth]{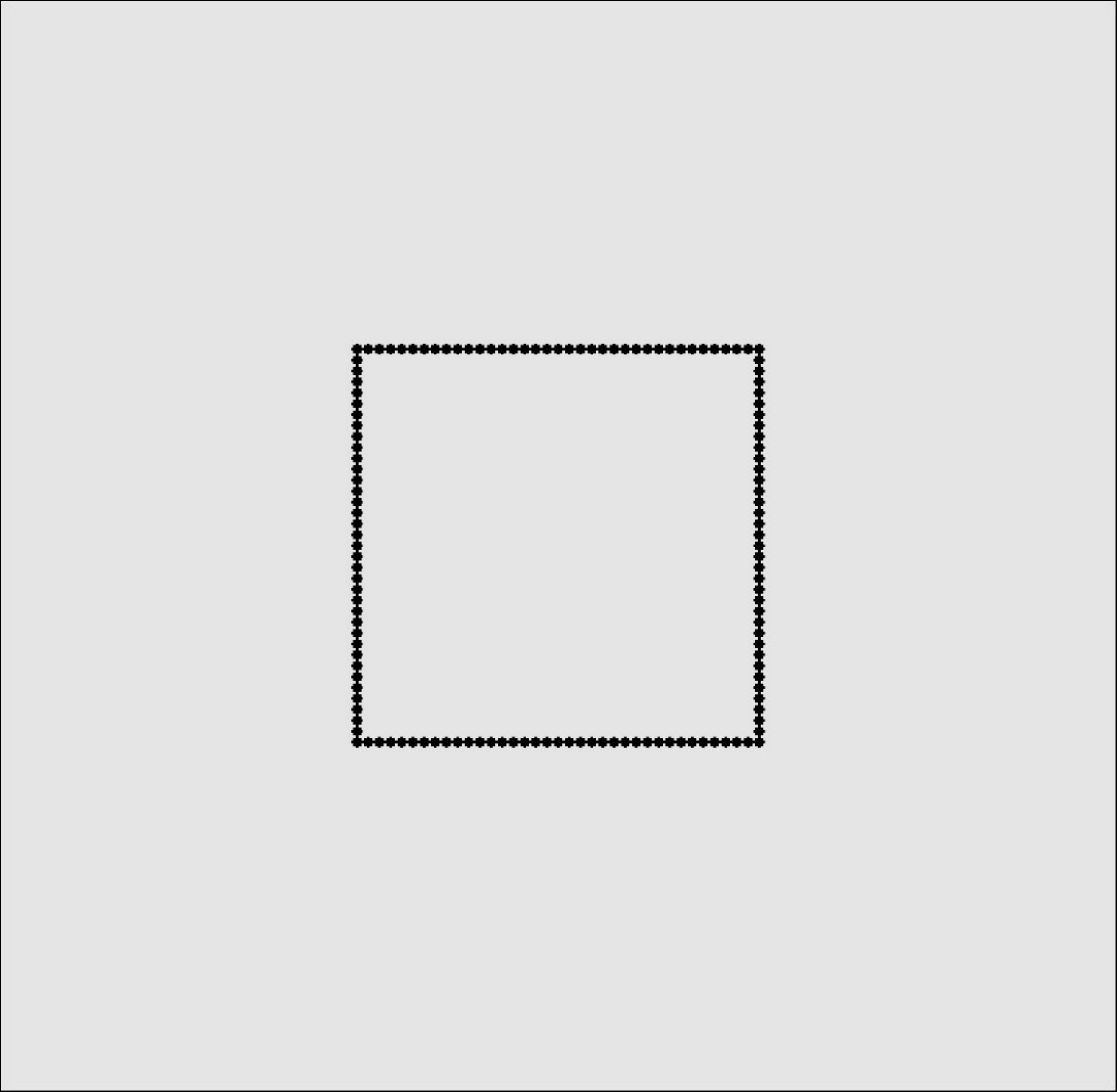} \qquad & 
  \includegraphics[width=.3\textwidth]{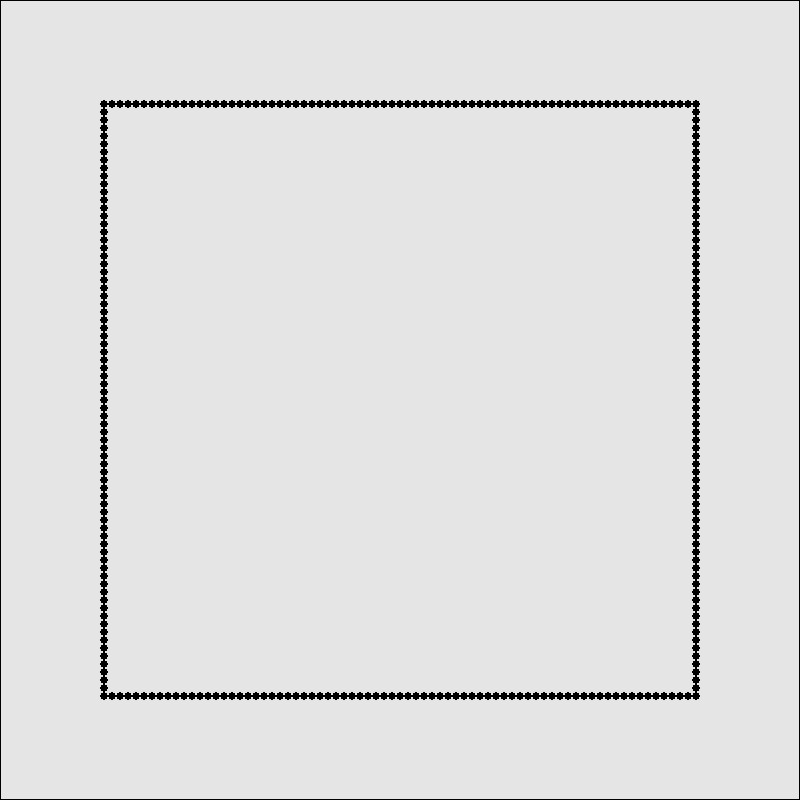} \qquad 
 \end{tabular}
\caption{FIM for HJB-1: active nodes at different steps}
\label{HJB-A}
\vskip-.2cm
\end{figure}

\begin{figure}[!h]
\centering
 \begin{tabular}{cccc}
  \includegraphics[width=.309\textwidth]{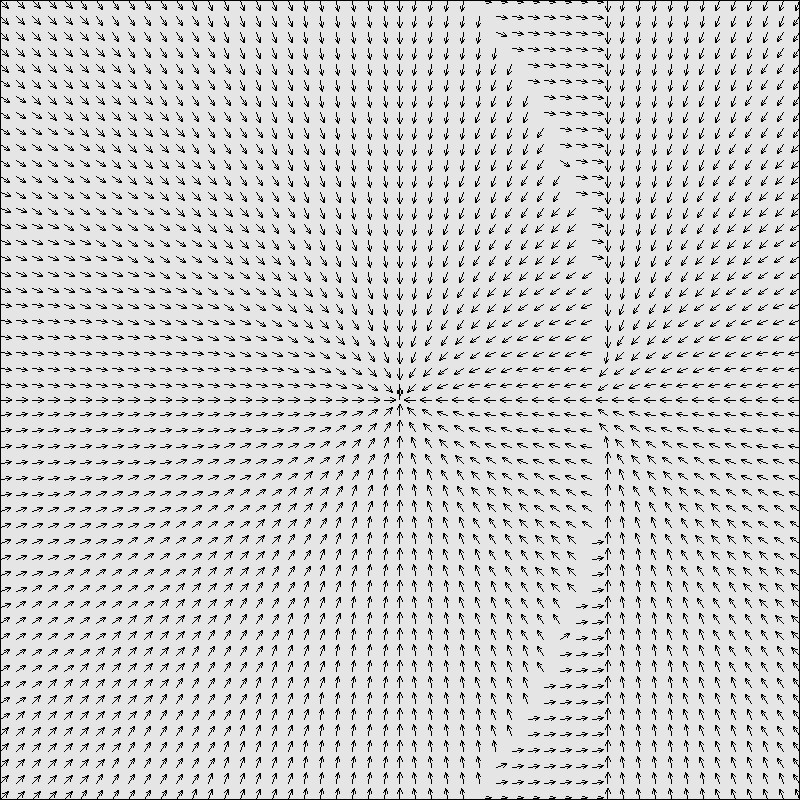}      \qquad & 
  \includegraphics[width=.3\textwidth]{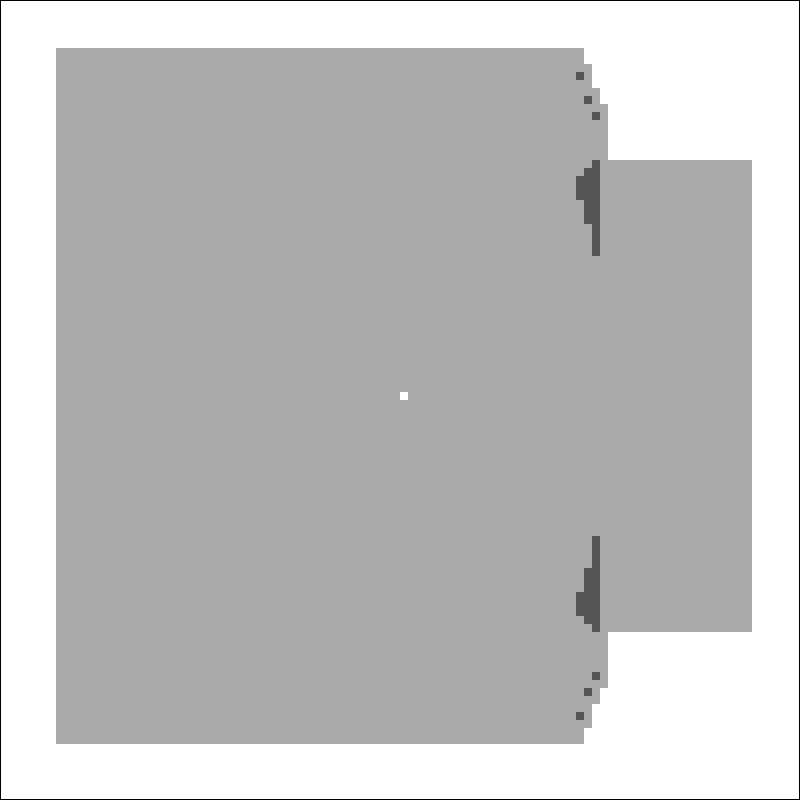} \qquad & 
  \includegraphics[width=.3\textwidth]{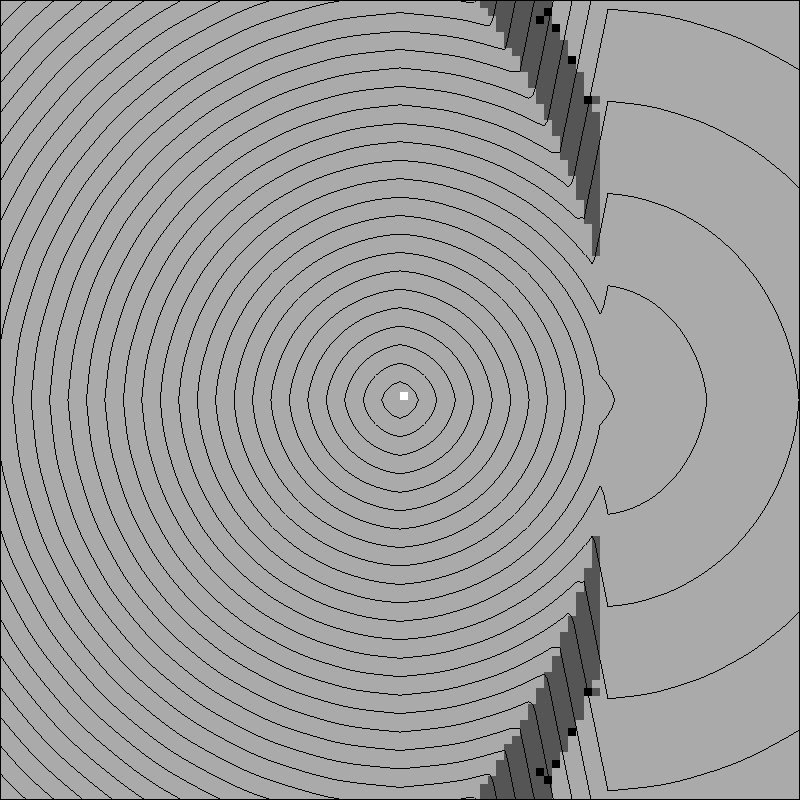} \qquad &
  \includegraphics[width=.025\textwidth,height=.3\textwidth]{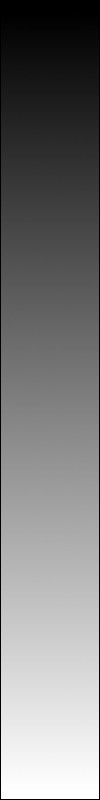}\begin{picture}(10,80)\put(2,0){0}\put(2,100){3}\end{picture}
 \end{tabular}
 \caption{FIM for HJB-2: Optimal controls (left), re-activation history at an intermediate step (center), re-activation history and level sets of the solution (right)}
 \label{HJB-B}
\vskip-.6cm
 \end{figure}
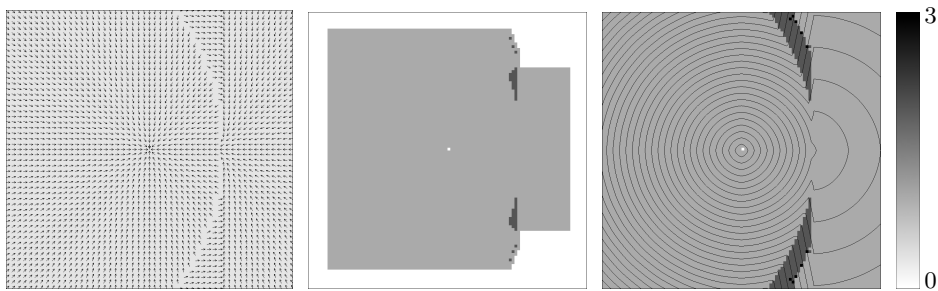

\begin{figure}[!h]
\centering
\begin{tabular}{cccc}
  \includegraphics[width=.309\textwidth]{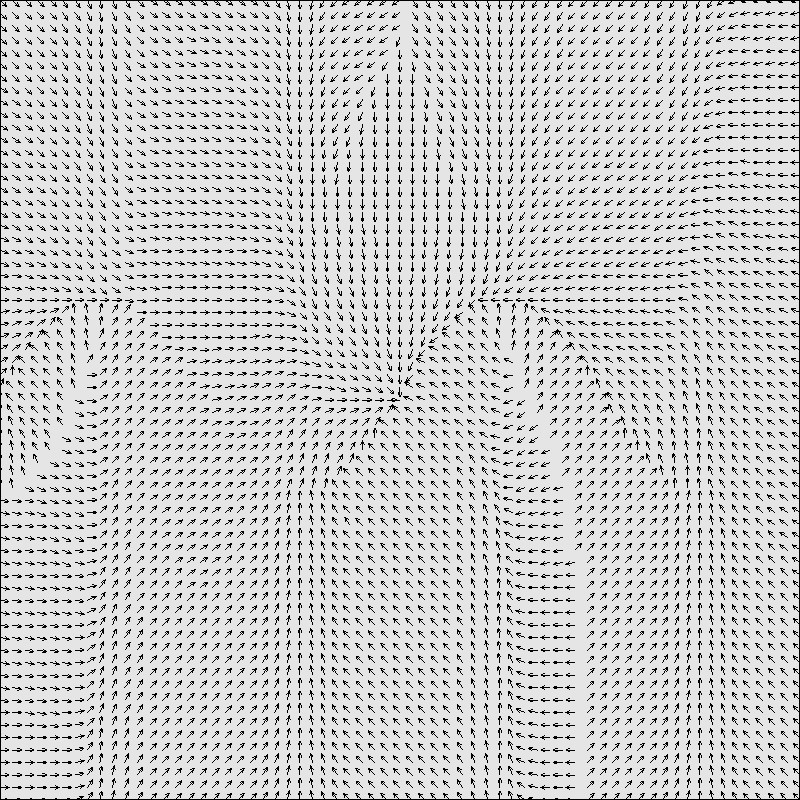} \qquad & 
  \includegraphics[width=.3\textwidth]{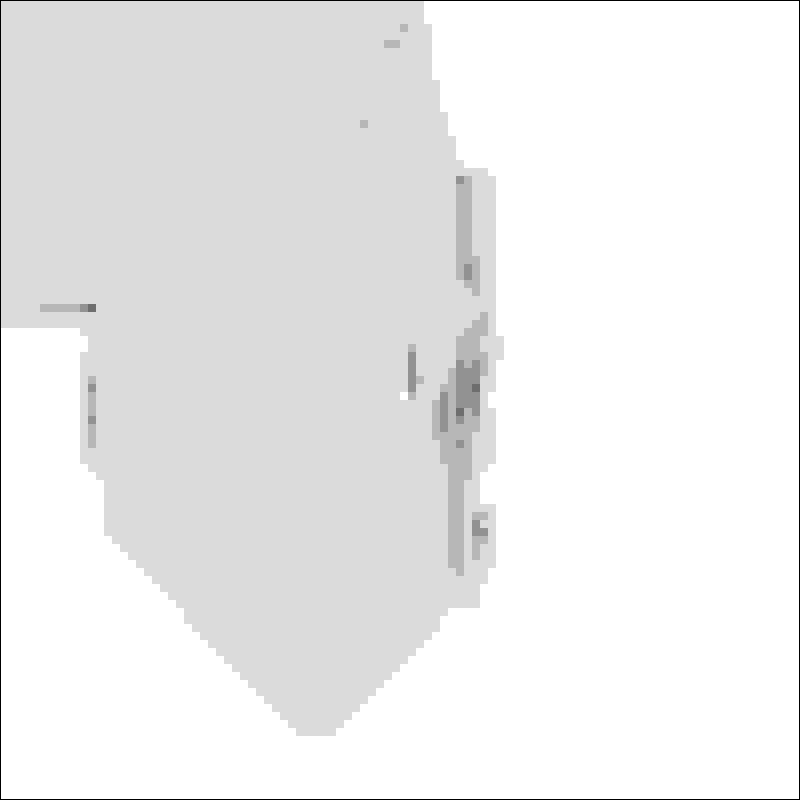} \qquad & 
  \includegraphics[width=.3\textwidth]{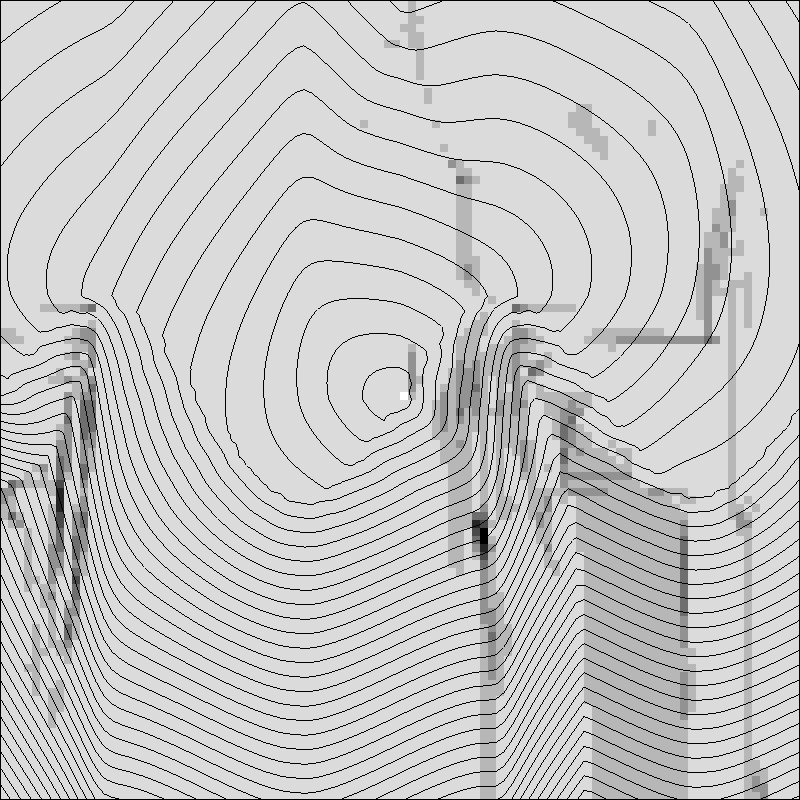} \qquad &
  \includegraphics[width=.025\textwidth,height=.3\textwidth]{figures/gradient.jpg}\begin{picture}(10,80)\put(2,0){0}\put(2,100){7}\end{picture}
  \end{tabular}
 \caption{FIM for HJB-4: Optimal controls (left), re-activation history at an intermediate step (center), re-activation history and level sets of the solution (right)}
 \label{seismic}
 \vskip-.7cm
\end{figure}

\begin{figure}[!h]
\centering
 \begin{tabular}{cccc}
  \includegraphics[width=.309\textwidth]{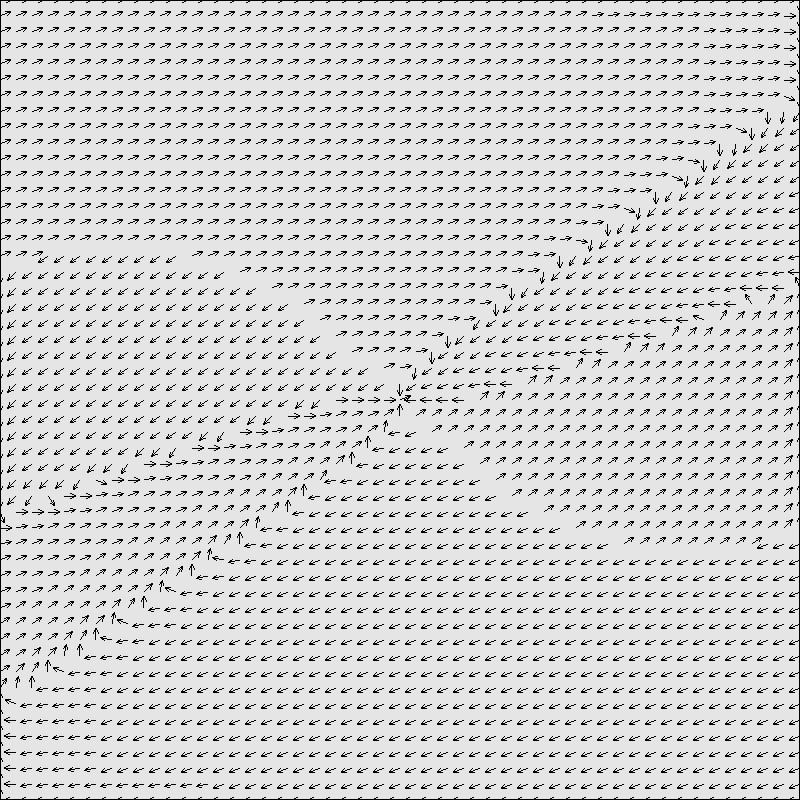} \qquad &     
   \includegraphics[width=.302\textwidth]{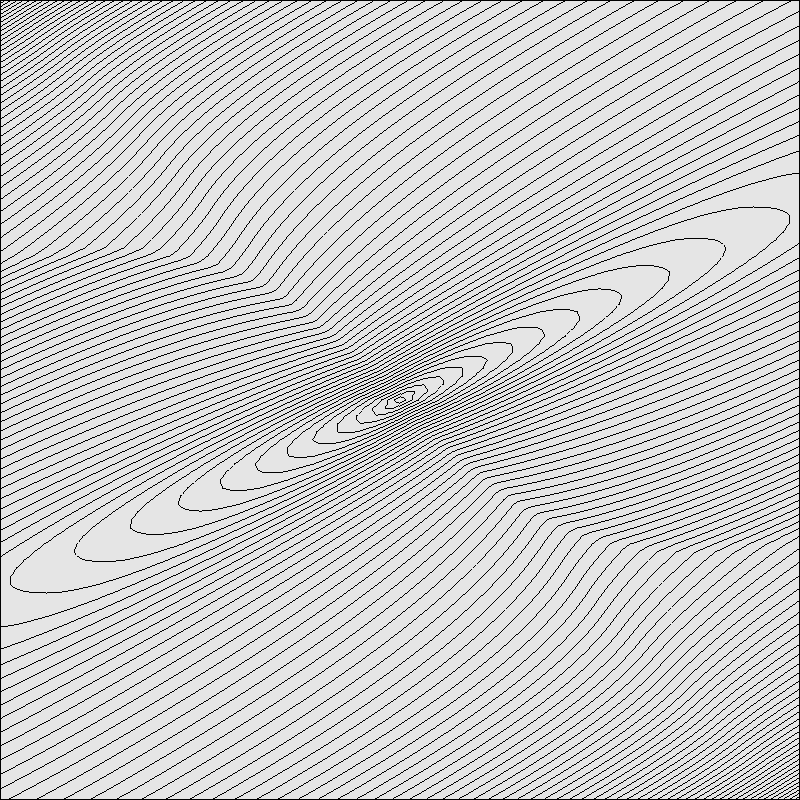} \qquad &    
  \includegraphics[width=.309\textwidth]{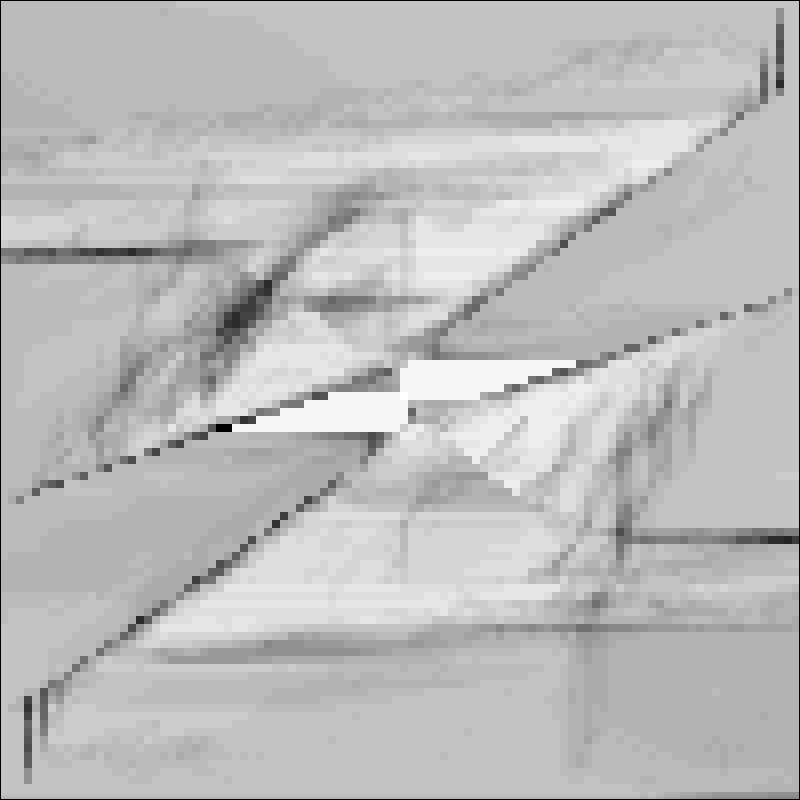} \qquad &
  \includegraphics[width=.025\textwidth,height=.3\textwidth]{figures/gradient.jpg}\begin{picture}(10,80)\put(2,0){0}\put(2,100){30}\end{picture}
  \end{tabular}
\caption{FIM for HJB-5: Optimal controls (left), level sets of the solution (center), re-activation history (right)}
\label{cravatta}
\end{figure}

In this section we compare the performance of FSM, FIM, UFSM{\footnotesize 1/4} and UFSM{\footnotesize 3/4} on the following equations:
\begin{center}
\begin{tabular}{|l|l|l|}
\hline Equation & Dynamics & Class
\\ \hline
\hline HJB-1 & $f(x,y,a)=a$ & ISO \& REG \\
\hline HJB-2 & $f(x,y,a)=(1+4\chi_{\{x>1\}})\,a$ & ISO \& $\neg$REG \\
\hline HJB-3 & $f(x,y,a)=m_{\lambda,\mu}(a)\,a$ & $\neg$ISO \& REG \\
\hline HJB-4 & $f(x,y,a)=F_2(x,y)m_{p(x,y),q(x,y)}(a)\,a$ & $\neg$ISO \& $\neg$REG \\
\hline HJB-5 & $f(x,y,a)=(1+|x+y|)m_{\lambda,\mu}(a)\,a$ & $\neg$ISO \& $\neg$REG \\
\hline
\end{tabular}
\end{center}
where we defined $\displaystyle m_{\lambda,\mu}(a)=(1+(\lambda\,a_1+\mu\,a_2)^2)^{-\frac{1}{2}}$ for $\lambda,\mu\in\R$ and we denoted by $\chi_S$ the characteristic function of a set $S$. 
Moreover, for $c_1,c_2,c_3,c_4>0$, we defined
$$
C(x)=c_1\sin\left(\frac{c_2\pi x}{c_3}+c_4\right)\,,\qquad
\big(F_1(x,y),F_2(x,y)\big)=\left\{\begin{array}{ll}
                                (0.5,1) & \mbox{if } y\leq C(x)\\
                                (2,3) & \mbox{otherwise}
                               \end{array}\right.\,,
$$
$$
M(x,y)=\sqrt{\frac{\frac{F_2^2(x,y)}{F_1^2(x,y)}-1}{1+{C^\prime}^2(x)}},\quad
p(x,y)=M(x,y)C^\prime(x),\quad q(x,y)=-M(x,y)\,.
$$
In all the following tests we set $\Omega=[-2,2]^2$ (except Test 4), the target $\mathcal{T}=(0,0)$ and the number of discrete controls $N_c=32$. 
Regarding FIM, we keep track of the history of the active list by counting the number $I_i$ of times the node $x_i$ enters the active list. 
The number $I_{\max}:=\max_i I_i$ gives a measure of the ``iterativeness'' of the method. 
\paragraph{Test 1. }
Here we solve equation HJB-1. Figure \ref{HJB-A} shows the evolution of FIM's active list. 
Differently from FMM, where the CONS set expands from the target following concentric circles (i.e. the level sets of the solution), 
here the active set moves following concentric squares (cf.\ the behavior of the CONS region of the \emph{safe method} studied in \cite{CCF14}). 
As one can expect $I_{\max}=1$, meaning that FIM behaves like a single-pass method.
Table \ref{tab:cpu-times} compares CPU times of the methods on different grids and the number of sweeps needed by sweeping methods to reach convergence. 
FSM converges in 4 sweeps for this equation, the additional sweep reported in Table \ref{tab:cpu-times} is the one required by the algorithm to check convergence. All the methods compute the same solution. In particular we see that FIM is slightly slower than FSM, as noted in \cite{JW08}. On the other hand, UFSM methods (both {\footnotesize 1/4} and {\footnotesize 3/4}) still converge in 5 sweeps, thus overcoming FSM.
\paragraph{Test 2. }
Here we solve equation HJB-2. Figure \ref{HJB-B} shows the optimal vector field $f(x,a^*)$ and the history of active nodes (in grey scale, where black corresponds to $I_{\max}$ and white to $0$). 
The maximal number of re-activation is $I_{\max}=3$ and re-activation of nodes appears for the first time close to the shock line, see Figure \ref{HJB-B}-center. 
This depends on the fact that  the active list is not an approximation of a level set of the solution and the equation HJB-2 falls in the $\neg$REG class.
Then FIM is not able to capture the shock properly (see Section \ref{sec:limitations}) in a single-pass fashion, but has to come back to recompute wrong values. 
Table \ref{tab:cpu-times} compares the methods. Results are similar to those of the previous test.
\vskip-.4cm
\paragraph{Test 3. }
Here we solve equation HJB-3 for $\lambda=10$ and $\mu=5$, namely the anisotropic eikonal equation, a well known example where FMM fails in computing the correct solution, due to the fact that characteristics do not coincide with gradient curves of the solution, see \cite{SV03}. 
In this case FIM let evolve the active list as for the eikonal equation (Test 1, Figure \ref{HJB-A}) and produces a maximal number of re-activation $I_{\max}=1$. Again, this means that equation HJB-3 can be successfully solved by a local single-pass method, as the \emph{safe method} introduced in \cite{CCF14} for the class REG. 
We refer to Table \ref{tab:cpu-times} for a comparison of the methods.
\paragraph{Test 4. }
Here we solve equation HJB-4 in $\Omega=[-0.5,0.5]^2$ for $c_1=0.1225$, $c_2=2$, $c_3=0.5$ and $c_4=0$, an example of class $\neg$ISO \& $\neg$REG coming from seismic imaging. 
It is a inhomogeneous anisotropic eikonal equation on a domain with two layers separated by a 
sinusoidal profile $C(x)$, with different constant anisotropy coefficients in each layer (given by the pairs $(F_1,F_2)=(0.5,1)$ and $(F_1,F_2)=(2,3)$). \\
All the methods compute the same solution, meaning that FIM can work for equations with substantial anisotropy and inhomogeneities (see also next test). 
Unexpectedly, also UFSM{\footnotesize 1/4} is able to correctly follow quite curved characteristics, see Figure \ref{seismic}-left).
Results in Table \ref{tab:cpu-times} show that sweeping methods need a large number of sweeps to reach convergence (even more for UFSMs, due to the control set reduction). 
This makes FIM be the fastest method.
The maximal number of re-activation for the active list is $I_{\max}=7$ and Figure \ref{seismic}-center/right shows that re-activation of nodes appears both close to the shocks and where the optimal field exhibits rapid changes of direction. 
\vskip-.3cm
\paragraph{Test 5. }
Here we solve HJB-5 for $\lambda=10$ and $\mu=5$. This is the hardest example of class $\neg$ISO \& $\neg$REG 
presented in \cite{CCF14}, where the shock (see the cubic-like curve in Figure \ref{cravatta}-left/center) and a strong anisotropy region meet at the target. 
Sweeping methods FSM and UFSM{\footnotesize 3/4} 
require much more sweeps with respect to the previous tests, while UFSM{\footnotesize 1/4} 
fails in computing the correct solution, confirming that the control set reduction to 1/4 of ball cannot be applied in any case.  

\begin{figure}[!h]
 \begin{tabular}{cccc}
  \includegraphics[width=.24\textwidth]{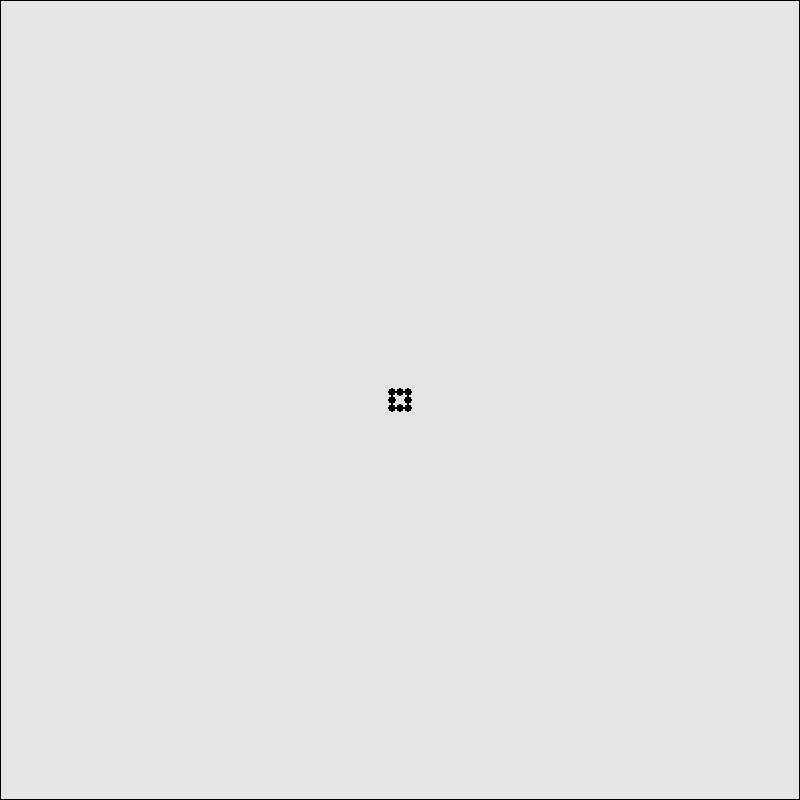} & 
  \includegraphics[width=.24\textwidth]{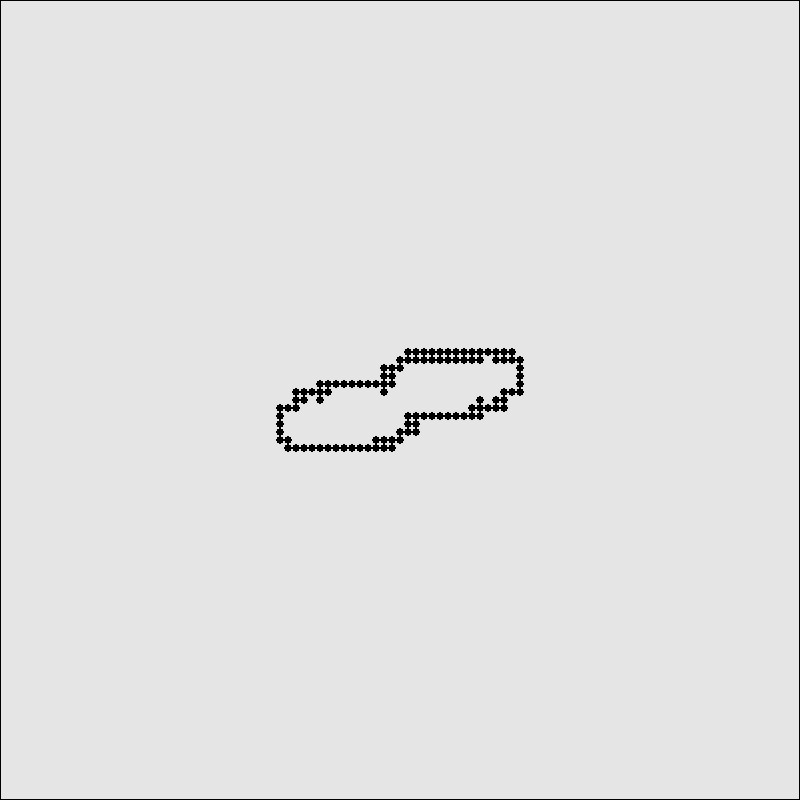} & 
  \includegraphics[width=.24\textwidth]{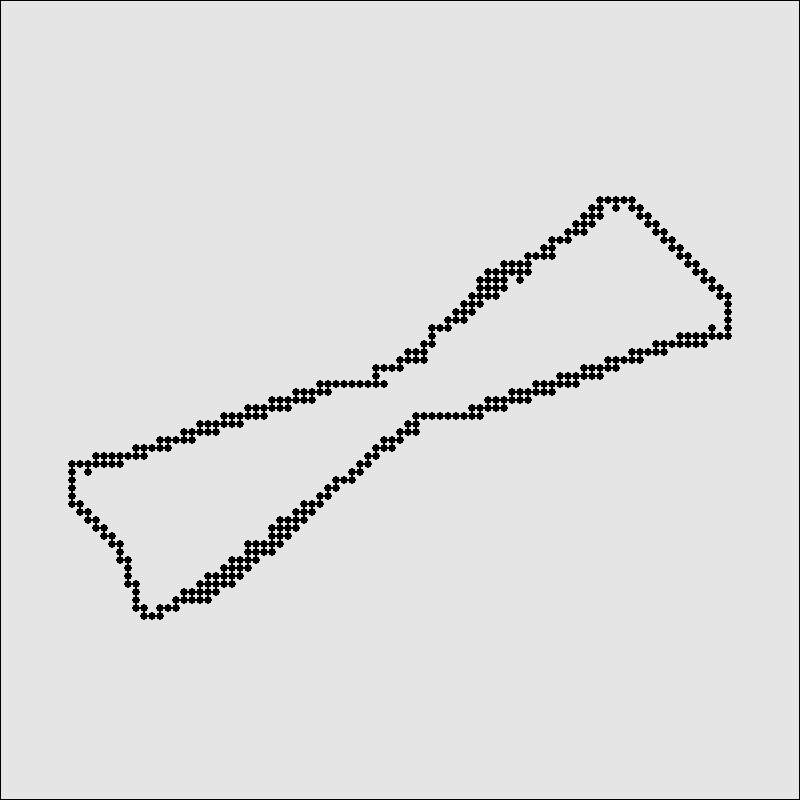} & 
  \includegraphics[width=.24\textwidth]{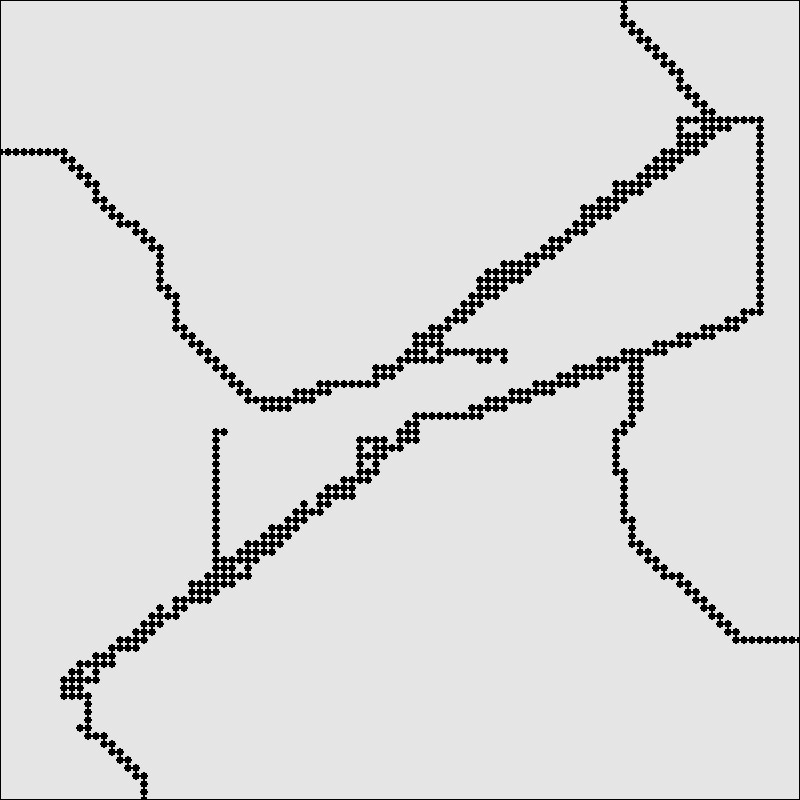} \\
  \includegraphics[width=.24\textwidth]{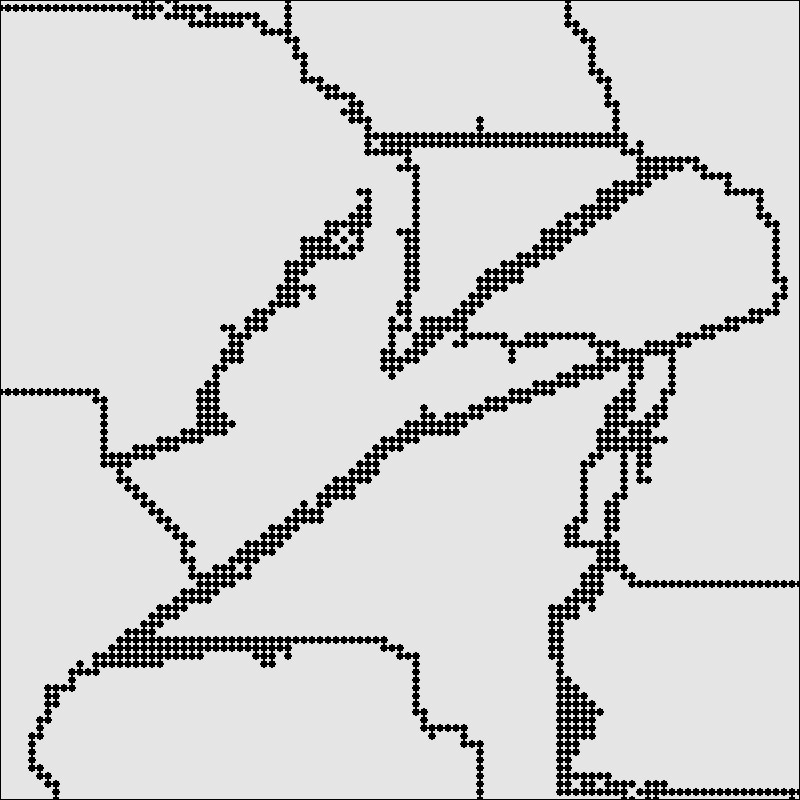} & 
  \includegraphics[width=.24\textwidth]{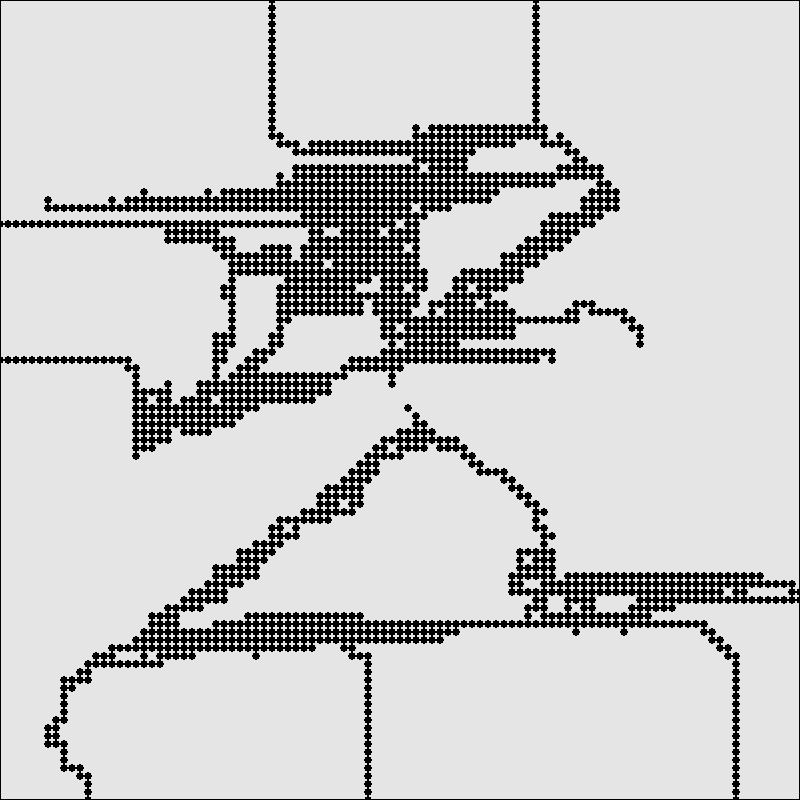} & 
  \includegraphics[width=.24\textwidth]{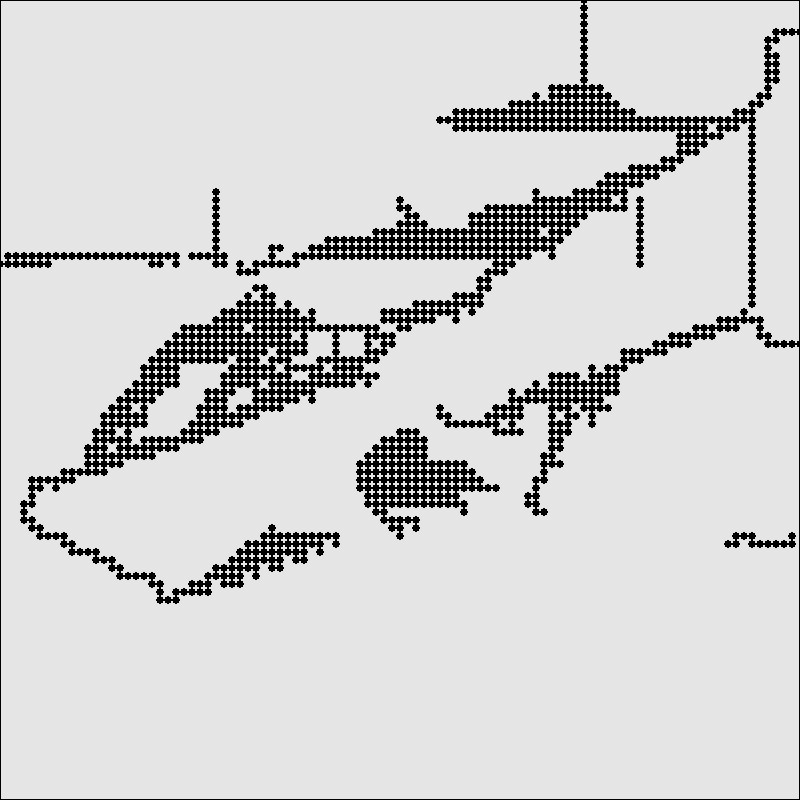} & 
  \includegraphics[width=.24\textwidth]{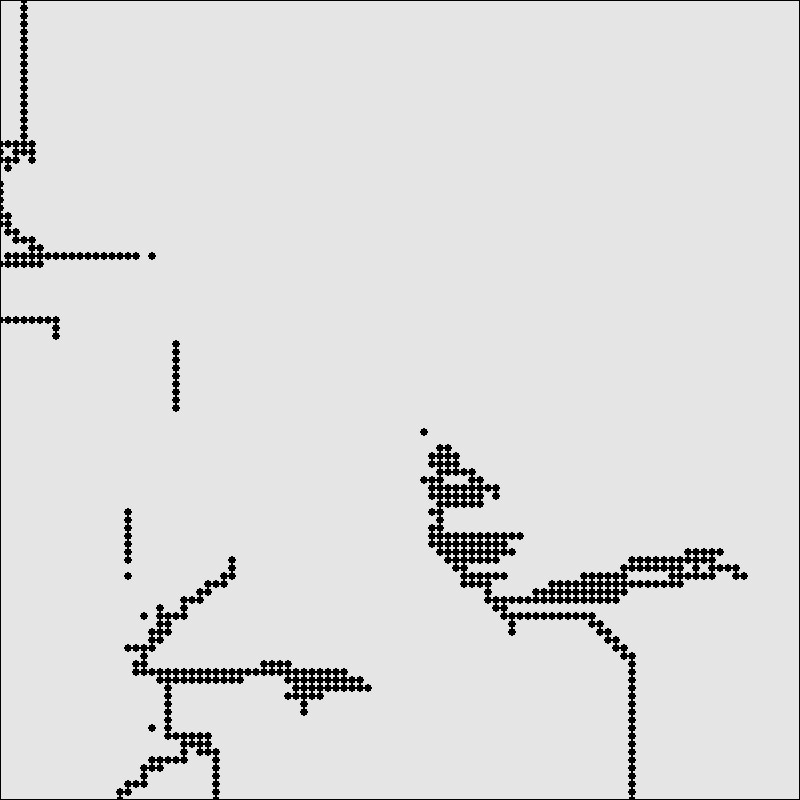}
 \end{tabular}
\caption{FIM for HJB-5: active nodes at different steps}\label{delirio}
\end{figure}

The maximal number of re-activation for FIM is $I_{\max}=30$ (see Figure \ref{cravatta}-right), whereas the evolution of the active list is extremely complicated and also produces \emph{regions of dimension two} (see Figure \ref{delirio}). Nevertheless, results in Table \ref{tab:cpu-times} shows that, as the grid increases, FIM is still the fastest method. 
The presence of a two-dimensional active list clearly proves that local single-pass methods cannot be applied since the enlargement of CONS is required (cf.\ the buffered fast marching method \cite{C09}).
\vskip-.6cm
\noindent {\bf Conclusions.} 
Tests performed in Section \ref{sec:tests} show that FIM can be successfully used for solving general HJB equations with no modifications 
with respect to the original algorithm. 
Moreover, FIM appears to be the fastest method in case of complicated $\neg$ISO \& $\neg$REG equations 
because of the large number of iterations needed by the sweeping methods.
Considering that the implementation of FIM is not harder than that of FSM, we think that, overall, FIM is the best method among the tested ones.

UFSM{\footnotesize 3/4} is always preferable to FSM, since the larger number of iterations needed for convergence are widely counterbalanced by the speedup 
for each single sweep. UFSM{\footnotesize 1/4} is instead not safely applicable for general equations.

Finally, the results of this paper confirm those of \cite{CCF14}. 
Complicated $\neg$ISO \& $\neg$REG equations require to pass through some nodes more than one time (cf., e.g., Figure \ref{cravatta}-right and \cite[Figure 7]{CCF14}), 
and exhibit two-dimensional regions in which every node depends on each other.

\begin{table}[h!]
\begin{center}
\begin{tabular}{|c|c|c|c|c|c|c|}
\hline
Equation & Grid & $\Delta x$ & FSM (sweeps) & ~FIM~ & UFSM{\footnotesize 1/4} (sweeps) & UFSM{\footnotesize 3/4} (sweeps)\\
\hline
 HJB-1 & 101 & 0.04  &   0.13 (5)  &    0.17   &    {\bf 0.04} (5)   &      0.10 (5) \\ 
 HJB-1 & 201 & 0.02  &   0.51 (5)  &    0.72   &    {\bf 0.15} (5)   &      0.40 (5) \\
 HJB-1 & 401 & 0.01  &   2.04 (5)  &    3.21   &    {\bf 0.58} (5)   &      1.51 (5) \\
\hline\hline
 HJB-2 & 101 & 0.04  &   0.16 (5) &     0.21   &    {\bf 0.04} (5)   &      0.12 (5) \\
 HJB-2 & 201 & 0.02  &   0.63 (5) &     0.87   &    {\bf 0.19} (5)   &      0.46 (5) \\
 HJB-2 & 401 & 0.01  &   2.46 (5) &     3.80   &    {\bf 0.70} (5)   &      1.84 (5) \\
\hline\hline
 HJB-3 & 101 & 0.04  &   0.31 (5) &     0.38   &    {\bf 0.09} (5)   &      0.23 (5) \\
 HJB-3 & 201 & 0.02  &   1.23 (5) &     1.56   &    {\bf 0.35} (5)   &      0.88 (5) \\
 HJB-3 & 401 & 0.01  &   4.88 (5) &     6.57   &    {\bf 1.38} (5)   &      3.53 (5) \\
\hline\hline
 HJB-4 & 101 & 0.01  &   5.72 (25) &    {\bf 1.93}   &    2.18 (34)  &      5.62 (34) \\
 HJB-4 & 201 & 0.005 &  22.70 (25) &    {\bf 7.68}   &    8.66 (34)  &     19.58 (30) \\
 HJB-4 & 401 & 0.0025 &  99.38 (28) &   {\bf 29.36}  &    34.07 (34) &     77.14 (30) \\ 
\hline\hline
 HJB-5 & 101 & 0.04  &   3.23 (53)  &   5.30    &   -   &     {\bf 2.38} (53)\\
 HJB-5 & 201 & 0.02  &  13.30 (55)  &   {\bf 3.32}    &   -   &     9.78 (55)\\
 HJB-5 & 401 & 0.01  &  52.93 (55)  &   {\bf 14.53}    &   -   &    39.16 (55)\\
\hline
\end{tabular}
\end{center}
\caption{CPU times (seconds) and number of sweeps. Fastest method in bold}\label{tab:cpu-times}
\vskip-.7cm
\end{table}

\end{document}